\documentstyle[12pt]{article}

\begin{document}
\newtheorem{guess}{Proposition }[section]
\newtheorem{theorem}[guess]{Theorem}
\newtheorem{lemma}[guess]{Lemma}
\newtheorem{corollary}[guess]{Corollary}
\def \ja {\vrule height 3mm width 3mm}

\vglue7pc
\centerline{\Large \bf On the topology of holomorphic bundles}
\vspace{.2 in}
\centerline{\Large Elizabeth Gasparim}
\centerline{\small Departamento de Matem\'atica, Universidade Federal de Pernambuco} 
\centerline {\small Cidade Universit\'aria, Recife, PE, BRASIL, 50670-901}
\centerline{\small gasparim@dmat.ufpe.br}

\vspace {5 mm}

\begin{abstract}

In this work we study 
the topology of holomorphic rank two bundles 
over complex surfaces.
We consider  bundles
that are constructed by glueing ``local'' holomorphic bundles
and we show that under certain conditions the 
topology of the bundle does not depend on the 
glueing. As a consequence we present a simple and new 
classification of bundles on blown-up surfaces.

\end{abstract}

\section{Introduction}

Let $X$ be a complex  manifold and let  
$A$ and $B$ be open sets that cover $X.$
Given holomorphic bundles 
 $E_A$ and $E_B$ defined
over the subsets $A$ and $B,$  we construct
holomorphic bundles over $X$
by glueing the bundles  $E_A$ and $E_B$ over $A \cap B.$
We  then compare the  topology of bundles 
given by different glueings.
 
Our main motivation is the study of bundles over blown-up 
surfaces. Therefore we will focus  our attention on
the cases where the intersection $A \cap B$ is 
biholomorphic to ${\bf C}^m - \{0\}.$
We then present  simple holomorphic and  topological
classifications of bundles on  blown-up surfaces.

\section{Glueing bundles over a  manifold}

In this section we
construct bundles over a complex manifold 
by glueing bundles defined on 
open subsets and then we
 compare their topology .\\

\noindent Consider the following data:
 
\vspace{ 2 mm}

\noindent i) a complex manifold
$X$ with $dim _{\bf C}X = m \ge 2$ 
and open subsets $A,B,$ and $C$ of $X$ satisfying
$X = A \cup B,\, A \cap B = C$ with $C$ biholomorphic 
to ${\bf C}^m - \{0 \}$

\vspace{2 mm}

\noindent ii) holomorphic vector bundles  $\pi_A: E_A \rightarrow A$ and
 $\pi_B : E_B \rightarrow B$
which are trivial when restricted to $C$

\vspace{2 mm}

\noindent iii) a trivialization $F: E_A \rightarrow
 C \times {\bf C}^n$ of $E_A|_C$ and trivializations 
 $G_i: E_B \rightarrow C \times {\bf C}^n, \, \, i = 0,1$ of $E_B|_C.$\\

Given the above data,
we define  bundles $E_i$  over $X$ by the formula
$$E_i = E_A \bigcup _{F = G_i} E_B 
 =  (E_A \bigsqcup  E_B)/ \sim$$
where for $x \in E_A$ and $y \in E_B$ we define
$x \sim y$ if $F(x) = G_i(y).$

\noindent We have the following result about the bundles $E_i.$

\begin{guess}\label{topequ} The topology of the bundle $E_i$ is independent 
of the glueing, that is $E_0$ and $E_1$ are topologically equivalent vector
bundles.
\end{guess}

\noindent To prove this we present some preliminary  lemmas.

\begin{lemma}\label{homtri}If the trivializations
 $G_0$ and $G_1$ are homotopic,
then  $E_0$ and $E_1$ are
topologically  equivalent vector bundles.
\end{lemma}

\noindent{\bf Proof}: 
By $G_0$ homotopic to $G_1$  we 
mean  that there is a one parameter family
of trivializations $G_t$ with $ t\in [0,1]$ 
taking $G_0$ to $G_1.$ 
 To prove the lemma  we consider the product  bundles 
$E_A \times I$ and 
$E_B \times I$ over  $A \times I$
and $B \times I$ respectively,
together with given   trivializations
$\widetilde{F}$ and $\widetilde{G}$
of $E_A \times I |_{C \times I}$
and  $E_B \times I |_{C \times I}$
defined by
$\widetilde{F}(a,t) = (F(a),t),$ for $(a,t) \in A \times I$
and $\widetilde{G}(b,t) = (G(b),t),$  for $(b,t) \in B \times I.$
We have that $\widetilde{G}(y,0) = G_0(y)$
and  $\widetilde{G}(y,1) = G_1(y).$
Defining  the bundle $E$ over $X \times I$ by
$E = (E_A \times I) \cup_{\widetilde{F} = \widetilde{G}} (E_B \times I)$
it immediately follows  that $E|_{X \times \{0\}} \simeq E_0$
and  $E|_{X \times \{1\}} \simeq E_1$
and consequently $E_0$ and $E_1$ are topologically equivalent. \hfill\ja

\begin{lemma}\label{tribun} Consider
 the trivial vector bundle $D = C \times {\bf C}^n$
over a complex space $C.$ Let  $G_0$ 
and $G_1$ be two trivializations of $D$ over $C$  and let
$\Phi:C \rightarrow GL(n,{\bf C})$ be 
the corresponding transition matrix, i.e. 
 $\Phi(c) G_0(c) = G_1(c).$ Then $G_0$ and $G_1$ are
homotopic if and only if $\Phi$ is nullhomotopic.
\end{lemma}

\noindent The proof is straightforward.

\begin{lemma}\label{nulhom} Any holomorphic map 
 $f: {\bf C}^m - \{0\} \rightarrow GL(n,{\bf C})$ 
for  $m\ge 2$
 is nullhomotopic.
\end{lemma}

\noindent {\bf Proof}: By Hartog's Theorem, $f$ extends to a holomorphic function 
$\widetilde{f}: {\bf C}^m \rightarrow M(n,{\bf C}),$ where $M(n,{\bf C})$
denotes the space of all $n \times n$ matrices with complex coefficients.
We claim that  $\widetilde{f}(0) \in GL(n,{\bf C}).$ 
In fact, if  $\widetilde{f}(0) \notin GL(n,{\bf C})$ then 
$det(\widetilde{f}(0)) = 0.$
But then $(det \circ \widetilde{f})^{-1} (0) = \{ 0\} \subset {\bf C}^m$
which is a contradiction, because  
$(det \circ \widetilde{f}): {\bf C}^m \rightarrow {\bf C}$ is holomorphic 
and the pre-image of a point by a
 holomorphic function in ${\bf C}^m$ is either empty 
or has codimension 1. Which proves the claim.

Thus, we can write
$\widetilde{f} : {\bf C}^m \rightarrow GL(n,{\bf C}).$  
Hence $f$ factors through ${\bf C}^m,$ 
i.e. $f = \widetilde{f} \circ i$ where $i$ is the inclusion
$i: {\bf C}^m -\{0\} \rightarrow {\bf C}^m.$ 
As ${\bf C}^m $ is contractible it follows that 
$f$ is nullhomotopic.

\vspace {3mm}
\noindent Alternative proof of Lemma \ref{nulhom}:
Consider the holomorphic function 
$g = det \circ f : {\bf C}^m - \{0 \} \rightarrow {\bf C}.$
By Hartog's theorem we have that both $g$ and $1/g$
extend to holomorphic functions defined on ${\bf C}^m.$
Let $\widetilde {g} $ be the extension of $g$ to 
${\bf C}^m. $ Then, because $1/g$
also extends as a holomorphic function to ${\bf C}^m,$
 it follows that  $\widetilde {g}(0) \neq 0.$ 
We have $g = i \circ \widetilde {g}$ where $i$ is 
the inclusion $i : {\bf C}^m - \{ 0\} \rightarrow {\bf C}^m$
and since ${\bf C}^m$ is contractible it follows that $g$
is nullhomotopic. \hfill\ja

\vspace{5 mm}

\noindent {\bf Proof of  Proposition \ref{topequ}}: Lemmas \ref{tribun}
 and \ref{nulhom}
imply that any two holomorphic trivializations
of the trivial bundle $({\bf C}^m - \{0 \}) \times {\bf C}^n$
over
 ${\bf C}^m - \{0 \}$ are homotopic and then
Lemma \ref{homtri}
  implies that the bundles obtained using these trivializations are 
topologically equivalent.\hfill\ja

\section{ Bundles on  Blown-up Surfaces } 
In this section we apply  Proposition \ref{topequ}
to give a simple topological description of bundles on some blown-up 
surfaces. Let us  consider the following case:

\vspace{2 mm}

\noindent  i) $X = \widetilde {S} $  is the blow-up 
 of a complex surface
$S$ at a point $P$  and $\ell$ is  the exceptional divisor 

\vspace{2 mm}

\noindent ii) the open subsets 
$A = N_{\ell}$ and $B = \widetilde {S} - \ell \simeq S - \{P\}$ 
are respectively a neighborhood
of the exceptional divisor $\ell$ and the complement of the 
exceptional divisor.

\vspace{2 mm}

\noindent iii) $A \cap B \simeq {\bf C}^2 - \{0\} .$\\

Some elementary
 examples of such surfaces are the blow up of the projective plane
${\bf P}^2$ at a point or the blow-up of a Hirzebruch surface 
$S_n = {\bf P}( {\cal O}(n) \oplus {\cal O})$ at a point.
Studying successive blow-ups on these basic surfaces leads
 to a similar topological classification of 
bundles on any rational surface. This just follows from the classification 
of rational surfaces, see Griffths and Harris [3].

To state our topological classification we
first quote some results from previous 
papers.

 We  write  $\widetilde {\bf C}^2  = U \cup V$,  where
$U =   {\bf C}^2 =\{(z,u)\},$\,\,
$V =   {\bf C}^2 = \{(\xi,v)\},$\,\,
$U \cap V = ({\bf  C} - \{0\})  \times   {\bf C}$
with the change of coordinates
$ (\xi,v) = (z^{-1 },zu).$

\begin{theorem}\label{matk=1}{\rm  [1, Thm. 2.1]}
  Let $E$ be a
 holomorphic rank two vector  bundle on
 $ \widetilde{\bf C}^2  $ with
vanishing first Chern class and
let $j$  be the non-negative integer that satisfies
$E_{\ell} \simeq {\cal O}(j) \oplus  {\cal O}(-j). $
Then  $E$  has a transition matrix
of the form

$$\left(\matrix {z^j & p \cr 0 &  z^{-j} \cr }\right)$$
from $U$ to $V,$  where 
$p$ is a polynomial given by
$$p = \sum_{i = 1}^{2j-2} \sum_{l = i-j+1}^{j-1}p_{il}z^lu^i.$$
\end{theorem}

\begin{theorem}\label{matgen}{\rm [2, Thm. 3.3]}
   Let $E$ be a
 holomorphic  rank two vector  bundle on
 ${\cal O}(-k)$ whose restriction to the zero section is  
$E_{\ell} \simeq  {\cal O}(j_1)  \oplus {\cal O}(j_2),$ with $j_1 \ge j_2.$   
Then  $E$  has a transition matrix
of the form
$$\left(\matrix {z^{j_1} & p  
                \cr   0 & z^{j_2}  
                \cr 
}\right)$$
from $U$ to $V,$  where the polynomial $p$ is given by

$$p = \sum_{i = 1}^{ \left[(j_1 - j_2 -2)/k\right]}
 \sum_{l = ki+j_2+1}^{j_1-1}p_{il}z^lu^i$$ and 
$p = 0$ if $j_1< j_2 +2.$
\end{theorem}

\begin{theorem}\label{triout}{\rm [2, Cor. 4.2]}
 Holomorphic bundles on the blow up
of a surface are trivial on a  neighborhood 
of the exceptional divisor minus the 
exceptional divisor.
\end{theorem}

As a consequence of Theorems \ref{matk=1} and \ref{triout} we have the 
following holomorphic classification of bundles on $\widetilde{S}.$

\begin{corollary}\label{holzer}
Every holomorphic rank two vector bundle over $\widetilde{S}$
with vanishing first Chern class 
is completely determined (up 
to isomorphism) by a 4-tuple $(E,j,p,\Phi)$ 
where $E$ is a holomorphic rank two bundle on $S$ with 
vanishing first Chern class, $j$ is a non-negative 
integer, $p$ is a polynomial, 
and $\Phi : {\bf C}^2 -\{0\} \rightarrow GL(2,{\bf C})$ is
a holomorphic map.
\end{corollary}

\noindent {\bf Proof}:
The essential ingredient here is that by theorem \ref{triout}
every holomorphic bundle 
on $\widetilde{S}$ is trivial on $N_\ell -\ell$ for some 
neighborhood $N_\ell$ of the exceptional divisor.
It follows that outside $\ell$ we may take a pull-back
bundle $\pi^*(E|_{S-p})$ of a holomorphic rank two bundle $E$
on $S$ with vanishing first Chern class
 and glue it to a bundle on $N_\ell$ using the function $\Phi.$
Now we use Theorem \ref{matk=1} to see that a bundle on $N_\ell$ 
is determined by a non-negative integer $j$ and a polynomial $p$ whose
form is explicitly known. \hfill\ja

\vspace{5 mm}

The corresponding classification for nonvanishing first 
Chern class is the following.

\begin{corollary}\label{holnon}
Every holomorphic rank two vector bundle over $\widetilde{S}$
is completely determined (up 
to isomorphism) by a 5-tuple $(E,j_1,j_2,p,\Phi)$ 
where $E$ is a holomorphic rank two bundle on $S,$
 $j_1$ and $j_2$ are 
integers, $p$ is a polynomial, 
and $\Phi : {\bf C}^2 -\{0\} \rightarrow GL(2,{\bf C})$ is
a holomorphic map.
\end{corollary}

\noindent The proof is analogous to the one for Corollary \ref{holzer}.

\section {Topology of bundles on $\widetilde{S}$}

We now deduce the topological counterparts of Corollaries 3.4 and 3.5.

\begin{corollary} \label{topzer}Every 
holomorphic rank two vector bundle over $\widetilde{S}$
with vanishing first Chern class 
is topologically  determined 
 by a triple $(E,j,p)$ 
where $E$ is a holomorphic rank two bundle on $S$ with 
vanishing first Chern class, $j$ is a non-negative
integer, and  $p$ is a polynomial.
\end{corollary}

\noindent {\bf Proof}: By Corollary \ref{holzer} we know that such a 
bundle is holomorphically
determined by a 4-tuple  $(E,j,p, \Phi)$ and 
Proposition \ref{topequ} shows that topologically the choice 
of the map
$\Phi$ is irrelevant.\hfill\ja

\vspace{5 mm}

A straightforward generalization of Corollary \ref{topzer}
 for the case of nonvanishing first Chern class 
is the following result.

\begin{corollary}\label{topnon} Every 
holomorphic rank two vector bundle over $\widetilde{S}$
is topologically  determined 
 by a 4-tuple $(E,j_1,j_2,p)$ 
where $E$ is a holomorphic rank two bundle on $S,$ 
 $j_1$ and $j_2$ are 
integers, and  $p$ is a polynomial.
\end{corollary}

\noindent The proof is analogous to the proof of \ref{topzer}.

\vspace {5 mm}

\noindent {\bf Examples}:
Let us write down some examples to clarify the statement of
Corollary \ref{topzer}.
First we fix a holomorphic rank two bundle $E$ 
with vanishing first Chern class over the 
surface $S.$
Then we would like to see which are the possible 
bundles $\widetilde {E}$ over $\widetilde {S}$ 
that are a pull-back of $E$ outside the exceptional divisor.
According to Corollary \ref{topzer}  
 such bundles are topologically given by a choice 
of an integer $j$ and a polynomial $p$ whose 
form is given  in Theorem \ref{matk=1} as
$p = \sum_{i = 1}^{2j-2} \sum_{l = i-j+1}^{j-1}p_{il}z^lu^i.$

If $j = 0,$ then $p = 0$ and it follows that
 $\widetilde {E} =\pi^* E$ is globally a pull-back bundle.
That is, applying Corollary \ref{holzer}.
we verify the well known fact that
bundles over a blown-up surface that are trivial 
when restricted to  the exceptional divisor 
are pull-backs. 

If $j = 1,$ then also $p = 0.$ In this case we
see that on a neighborhood $N_\ell$  of $\ell$ the
bundle is ${\cal O}(1) \oplus {\cal O}(-1).$
If follows that  all holomorphic  bundles $\widetilde {E}$
over $\widetilde {S}$ whose 
restriction to the exceptional is ${\cal O}(1) \oplus {\cal O}(-1).$
 are topologically equivalent. However, clearly these bundles 
are not pull-backs and are not topologically equivalent 
to any of the  bundles we obtained in the previous case.

If $j = 2,$ then $p = (p_{10} + p_{11}z)u + p_{21} zu^2$ depends on three 
complex parameters. However these are not effective parameters 
in the sense that some different choices of
the polynomial $p$ will give isomorphic bundles over 
$N_\ell$ (holomorphically and hence also topologically)
and therefore will also lead to globally isomorphic bundles 
over $\widetilde {S}.$
This leads us immediately to the question of
 determining the ``local moduli space''
structure. That is, to see for a fixed value of $j$
what polynomials define isomorphic bundles over 
$N_\ell.$
We call ${\cal M}_j$ the moduli space of
isomorphism classes of  bundles 
on $N_\ell$ whose restriction to $\ell$
equals  ${\cal O}(j) \oplus {\cal O}(-j).$
 The answer to the local moduli question
is given by the following results.

\begin{theorem}\label{modj=2}{\rm [1, Thm. 3.4]}
 The moduli space ${\cal M}_2$ is homeomorphic to 
 the union ${\bf P}^1 \cup \{q_1,q_2\},$ 
of a complex projective plane ${\bf P}^1$ and two points,
  with  a basis  of open sets given by
$${\cal U}  \cup \{q_1,U : U \in {\cal U} - \phi \} \cup
\{q_1,q_2,U  : U \in {\cal U} - \phi \}  $$
where ${\cal U}$ is a basis for the standard topology of  ${\bf P}^1.$
\end{theorem}

\begin{theorem}\label{modgen}{\rm [1, Thm.3.5]}
 The generic set of the moduli space ${\cal M}_j$ is a complex 
projective space of dimension $2j-3.$
\end{theorem}

\noindent{ \bf Examples}:
Let us continue the analysis of the case $j = 2$ started 
in the preceding example.
We have seen that for $j = 2$ the polynomial $p$
is given by three complex parameters. However,
by Theorem  \ref{modj=2} we see that nonequivalent choices of $p$
are parametrized by a 
non-Hausdorff space formed by a  projective line
${\bf P}^1$   
with two extra points.
Therefore, it follows from Corollary \ref{holzer}
that for each fixed choice of glueing $\Phi,$
isomorphism classes of  bundles are parametrized by
 ${\bf P}^1$ (with the standard topology)
plus two extra points.
However it is simple to see that 
any such choices will produce topologically equivalent bundles.

For each value of the integer $j$ we 
can reproduce an  analysis similar to the ones 
in  the previous examples.
For a chosen bundle $E$ over $S$ and a particular choice
of glueing $\Phi$ we have generically a 
projective space $P^{2j-3}$ parametrizing 
nonisomorphic bundles $\widetilde {E}.$
Details for the topology as well as
explicit calculations of  Chern classes 
for these bundles will appear in a 
subsequent paper.

Let us represent holomorphic bundles
on the blown-up surface $\widetilde{S}$ by
$(E,j,p,\Phi),$ according to Corollary \ref{holzer}.
Then we have just proved the following.

\begin{theorem} Let $E$ be a holomorphic rank two 
bundle with vanishing 
first Chern class.
For fixed $\Phi$
and $j$ the family of isomorphism
classes of holomorphic bundles on 
 $\widetilde {S}$ which are a pull-back 
of $E$ outside the exceptional divisor
is generically parametrized by ${\bf P}^{2j-3}.$
\end{theorem}

\centerline{\bf Acknowledgments}
I am happy to thank Pedro Ontaneda for all the help.

\end{document}